\documentclass[11pt]{article}
\usepackage{latexsym,amsfonts,amssymb,amsmath,amsthm}
\usepackage{graphicx}

\usepackage[usenames,dvipsnames]{color}
\usepackage{ulem}

\begin{document}
\setlength{\baselineskip}{16pt}

\newtheorem{theorem}{Theorem}[section]
\newtheorem{lemma}{Lemma}[section]
\newtheorem{proposition}{Proposition}[section]
\newtheorem{definition}{Definition}[section]
\newtheorem{example}{Example}[section]
\newtheorem{corollary}{Corollary}[section]

\newtheorem{remark}{Remark}[section]

\numberwithin{equation}{section}

\def\p{\partial}
\def\I{\textit}
\def\R{\mathbb R}
\def\C{\mathbb C}
\def\u{\underline}
\def\l{\lambda}
\def\a{\alpha}
\def\O{\Omega}
\def\e{\epsilon}
\def\ls{\lambda^*}
\def\D{\displaystyle}
\def\wyx{ \frac{w(y,t)}{w(x,t)}}
\def\imp{\Rightarrow}
\def\tE{\tilde E}
\def\tX{\tilde X}
\def\tH{\tilde H}
\def\tu{\tilde u}
\def\d{\mathcal D}
\def\aa{\mathcal A}
\def\DH{\mathcal D(\tH)}
\def\bE{\bar E}
\def\bH{\bar H}
\def\M{\mathcal M}
\renewcommand{\labelenumi}{(\arabic{enumi})}

\def\disp{\displaystyle}
\def\undertex#1{$\underline{\hbox{#1}}$}
\def\card{\mathop{\hbox{card}}}
\def\sgn{\mathop{\hbox{sgn}}}
\def\exp{\mathop{\hbox{exp}}}
\def\OFP{(\Omega,{\cal F},\PP)}
\newcommand\JM{Mierczy\'nski}
\newcommand\RR{\ensuremath{\mathbb{R}}}
\newcommand\CC{\ensuremath{\mathbb{C}}}
\newcommand\QQ{\ensuremath{\mathbb{Q}}}
\newcommand\ZZ{\ensuremath{\mathbb{Z}}}
\newcommand\NN{\ensuremath{\mathbb{N}}}
\newcommand\PP{\ensuremath{\mathbb{P}}}
\newcommand\abs[1]{\ensuremath{\lvert#1\rvert}}

\newcommand\normf[1]{\ensuremath{\lVert#1\rVert_{f}}}
\newcommand\normfRb[1]{\ensuremath{\lVert#1\rVert_{f,R_b}}}
\newcommand\normfRbone[1]{\ensuremath{\lVert#1\rVert_{f, R_{b_1}}}}
\newcommand\normfRbtwo[1]{\ensuremath{\lVert#1\rVert_{f,R_{b_2}}}}
\newcommand\normtwo[1]{\ensuremath{\lVert#1\rVert_{2}}}
\newcommand\norminfty[1]{\ensuremath{\lVert#1\rVert_{\infty}}}
\newcommand{\ds}{\displaystyle}

\title{Criteria for the Existence of Principal Eigenvalues of Time Periodic Cooperative Linear Systems with Nonlocal Dispersal}

\author{Xiongxiong Bao\\
School of Science,\\
Chang'an University,\\
Xi'an, Shaanxi 710064, China\\
\\
   Wenxian Shen\\
   Department of Mathematics and Statistics\\
Auburn University\\
Auburn University, AL 36849\\
\\
E-mail address:
baoxx09@lzu.edu.cn,\, wenxish@auburn.edu }

\date{}

\maketitle

\vspace{-0.1in}

\noindent {\bf Abstract.} The current paper {{establishes}} criteria for the existence of principal eigenvalues of time periodic  cooperative linear nonlocal dispersal    systems with Direchlet type, Neumann type or periodic type boundary conditions. It is  shown that such a nonlocal dispersal system has a principal eigenvalue in the following cases: the nonlocal dispersal distance is sufficiently small; the spatial inhomogeneity satisfies a so called vanishing condition;  or the spatial inhomogeneity is nearly globally homogeneous. Moreover, it is shown that the principal eigenvalue of a time periodic cooperative linear nonlocal dispersal  system (if it exists) is algebraically simple. A linear nonlocal dispersal  system may not have a principal eigenvalue. The results established in the current paper extend those in literature for time independent or periodic nonlocal dispersal equations
to time periodic  cooperative nonlocal dispersal systems and {{will}} serve as {{a}} basic tool for the study of cooperative nonlinear systems with nonlocal dispersal.

\medskip


\medskip


\vspace{-0.1in}
\section{Introduction}

\vspace{-0.05in}

{{Internal dispersal (diffusion), which describes the movements
or  interactions of the organisms in the underlying systems,
 occurs in many  evolution systems arising in applied sciences.
There are several disparate approaches to model internal dispersal.
In case that the movements or  interactions of the organisms occurs randomly
between adjacent spatial locations,  differential operators  such as $\Delta$, have  been
widely used to model the dispersal,
see, for example,  \cite{ArWe75}, \cite{ArWe78},  \cite{Cantrell2003}, \cite{Hess1991},   \cite{Shigesada1997},
\cite{Smith1995}, etc..  They are frequently referred to
as {\it random dispersal operators}.
In contrast, integral operators, which are referred to as {\it nonlocal dispersal operators}, are often used to describe the movements or
interactions of the organisms  in an underlying system which occur between adjacent as well as non-adjacent spatial locations, see, for example,
\cite{Fi03}, \cite{Grinfeld2005},  \cite{Huston2003}, \cite{Huston2006}, etc.. For the study of various dynamical aspects of nonlocal dispersal evolution equations, the reader is referred to
 \cite{Bates2007},  \cite{Coville2006NA}-\cite{Coville2008}, \cite{Coville2013A},
\cite{Grinfeld2005}, \cite{Hetzer2012}, \cite{Huston2003}-\cite{Pan2010},  \cite{Rawal2014}, \cite{ShSh1}-\cite{ShSh4},
\cite{Shen20111}, \cite{Shen2012}.

It is of great importance  for
the study of nonlinear dispersal evolution equations/systems
to establish a spectral theory for
linear  dispersal evolution equations/systems.
The objective of the current paper is to  develop the  principal eigenvalue theory for
the following three  linear cooperative systems with nonlocal
dispersal and time periodic dependence:}}

\vspace{-0.1in}\begin{equation}
\label{DC-eq} \mathbf{u}_t=\int_{D}\kappa
(y-x)\mathbf{u}(t,y)dy-\mathbf{u}(t,x)+A_{1}(t,x)\mathbf{u}(t,x),\quad
x\in\bar D,
\vspace{-0.05in}\end{equation}
where $D\subset \Bbb{R}^{N}$ is a smooth bounded domain,
$\mathbf{u}(t,x)\in \Bbb{R}^{K}\ (K>1)$ and
$A_{1}(t,x)=(a_{kj}(t,x))_{K\times K}$ with $a_{kj}(t,x)$ being continuous in $t\in\RR$ and $x\in\bar D$ and
$a_{kj}(t+T,x)=a_{kj}(t,x)$ for $k,j=1,...,K$,
\vspace{-0.1in}\begin{equation}
\label{NC-eq} \mathbf{u}_t=\int_{D}\kappa
(y-x)[\mathbf{u}(t,y)-\mathbf{u}(t,x)]dy+A_{2}(t,x)\mathbf{u}(t,x),\quad
x\in\bar D,
\vspace{-0.1in}\end{equation}
where $ D\in\Bbb{R}^{N}$ and $A_{2}(t,x)$ {{has}} the same properties as $A_1(t,x)$,
and
\vspace{-0.05in}\begin{equation}
\label{periodic-eq} \mathbf{u}_t=\int_{\RR^N}\kappa
(y-x)[\mathbf{u}(t,y)-\mathbf{u}(t,x)]dy+A_{3}(t,x)\mathbf{u}(t,x),\quad
x\in\RR^N,
\vspace{-0.1in}\end{equation}
where $A_{3}(t,x)=(a_{kj}(t,x))_{K\times K}$ with $a_{kj}(t,x)$ being continuous in $t\in\RR$ and $x\in\RR^N$ and
$a_{kj}(t+T,x)=a_{kj}(t,x+p_{l}\mathbf{e}_{l})=a_{kj}(t,x)$ for
$p_{l}>0$, $\mathbf{e}_{l}=(\delta_{l1},...,\delta_{lN})$ with
$\delta_{ln}=1$ if $l=n$ and $\delta_{ln}=0$ if $l\neq n$,
$l,n=1,...,N$; $\kappa(\cdot)$ in \eqref{DC-eq}-\eqref{periodic-eq}
is a nonnegative $C^{1}$ function with compact support,
$\kappa(0)>0$, $\int_{\RR^N}\kappa(z)dz=1$.

{{In order to}} state the standing assumptions,  let
\vspace{-0.07in}\begin{equation}
\label{d-eq}
D_1=D_2=\bar D,\quad D_3=\RR^N.
\vspace{-0.07in}\end{equation}
Throughout this paper, we assume that

\begin{description}
\vspace{-0.1in}\item[(H1)] For each given $1\leq i\leq 3$, $A_{i}(\cdot,\cdot)$ is  cooperative in the sense
that for any $t\in\RR$, $x\in D_i$, and $1\leq k\neq j\leq K$,
$a_{kj}(t,x)\geq 0$.

\item[(H2)] For each given $1\leq i\leq 3$, $A_{i}(\cdot,\cdot)$ is strongly
irreducible in the sense that  for
any two nonempty subsets $S$, $S^{\prime}$ of $\{1,2,...,K\}$ which form a
partition of $\{1,2,...,K\}$ and any $t\in\RR$ and $x\in D_i$  there exist $k\in
S$, $j\in S^{\prime}$ satisfying that $|a_{kj}(t,x)|>0$.
\end{description}

\vspace{-0.1in} The  following eigenvalue problems are associated {{with the}} coupled nonlocal dispersal
systems \eqref{DC-eq}, \eqref{NC-eq}, and \eqref{periodic-eq}, respectively,
\begin{equation}\label{DC-eigen}
\begin{cases}
 -\mathbf{u}_t+\int_{D}\kappa
(y-x)\mathbf{u}(t,y)dy-\mathbf{u}(t,x)+A_{1}(t,x)\mathbf{u}(t,x)=\lambda\mathbf{u}(t,x),\quad
x\in\bar D\\
\mathbf{u}(t+T,x)=\mathbf{u}(t,x),
\end{cases}
\end{equation}
\begin{equation}
\label{NC-eigen}
\begin{cases}
 -\mathbf{u}_t+\int_{D}\kappa
(y-x)[\mathbf{u}(t,y)-\mathbf{u}(t,x)]dy+A_{2}(t,x)\mathbf{u}(t,x)=\lambda\mathbf{u}(t,x),\quad
x\in\bar D\\
\mathbf{u}(t+T,x)=\mathbf{u}(t,x),
\end{cases}
\end{equation}
and
\begin{equation}
\label{periodic-eigen}
\begin{cases}
-\mathbf{u}_t+\int_{\RR^N}\kappa
(y-x)[\mathbf{u}(t,y)-\mathbf{u}(t,x)]dy+A_{3}(t,x)\mathbf{u}(t,x)=\lambda\mathbf{u}(t,x),\quad
x\in\RR^N\\
\mathbf{u}(t+T,x)=\mathbf{u}(t,x+p_l{\bf e_l})=\mathbf{u}(t,x),\,\,\, l=1,2,\cdots,N.
\end{cases}
\end{equation}

The so-called principal spectrum point of \eqref{DC-eigen}  (resp. \eqref{NC-eigen},
\eqref{periodic-eigen}) plays a special role in the study of the related nonlinear systems.
Roughly, let $\sigma(A_1)$ (resp. $\sigma(A_2)$, $\sigma(A_3)$) be the spectrum of the eigenvalue problem  \eqref{DC-eigen}  (resp. \eqref{NC-eigen},
\eqref{periodic-eigen}). Let $\lambda_i(A_i)=\sup\{{\rm Re}\lambda\,|\,\lambda\in \sigma(A_i)\}$.
 $\lambda_{1}(A_{1})$ (resp. $\lambda_{2}(A_{2})$,
$\lambda_{3}(A_{3})$) is called the {\it principal spectrum point} of the
eigenvalue problem \eqref{DC-eigen} (resp. \eqref{NC-eigen},
\eqref{periodic-eigen}) (see Definition \ref{def2.2} for details). $\lambda_{1}(A_{1})$
(resp. $\lambda_{2}(A_{2})$, $\lambda_{3}(A_{3})$) is called the
{\it principal eigenvalue} of \eqref{DC-eigen} (resp. \eqref{NC-eigen},
\eqref{periodic-eigen}) if it is an isolated eigenvalue with a
positive eigenfunction  (see Definition \ref{def2.2} for details).

Observe that  \eqref{DC-eigen}, \eqref{NC-eigen}, and
\eqref{periodic-eigen} can be viewed as the nonlocal dispersal counterparts of the following eigenvalue problems with random dispersal and Dirichlet,
Neumann, and periodic boundary conditions, respectively,
\begin{equation}\label{DC-random-eigen}
\begin{cases}
 -\mathbf{u}_t+\Delta \mathbf{u}(t,x)+A_{1}(t,x)\mathbf{u}(t,x)=\lambda\mathbf{u}(t,x),\quad
x\in D\\
\mathbf{u}(t,x)=0,\quad x\in \p D\\
\mathbf{u}(t+T,x)=\mathbf{u}(t,x),
\end{cases}
\end{equation}
\begin{equation}
\label{NC-random-eigen}
\begin{cases}
 -\mathbf{u}_t+\Delta \mathbf{u}(t,x)+A_{2}(t,x)\mathbf{u}(t,x)=\lambda\mathbf{u}(t,x),\quad
x\in D\\
\frac{\p \mathbf{u}}{\p n}=0,\quad x\in\p D\\
\mathbf{u}(t+T,x)=\mathbf{u}(t,x),
\end{cases}
\end{equation}
and
\begin{equation}
\label{periodic-random-eigen}
\begin{cases}
-\mathbf{u}_t+\Delta \mathbf{u}(t,x)+A_{3}(t,x)\mathbf{u}(t,x)=\lambda\mathbf{u}(t,x),\quad
x\in\RR^N\\
\mathbf{u}(t+T,x)=\mathbf{u}(t,x+p_l{\bf e_l})=\mathbf{u}(t,x),\,\,\, l=1,2,\cdots,N.
\end{cases}
\end{equation}
In fact, it is proved for $K=1$ that the principal eigenvalue of  \eqref{DC-random-eigen}, \eqref{NC-random-eigen}, and
 \eqref{periodic-random-eigen} can be approximated by the principal spectrum
points of \eqref{DC-eigen}, \eqref{NC-eigen},  and
\eqref{periodic-eigen}, respectively, by properly rescaling the kernels  and that the
initial value problems of random dispersal evolution equations can also be
approximated by the initial value problems of the related nonlocal dispersal evolution equations
with properly rescaled kernels  (see
\cite{Cortazar2009}, \cite{Cortazar2008} and
\cite{ShenDCDS2015}). Hence we may say that \eqref{DC-eq},
\eqref{NC-eq}, and \eqref{periodic-eq} correspond to problems having Dirichlet type
boundary conditions, Neumann type boundary conditions,  and periodic
boundary conditions, respectively.


When $K=1$,   there are many results on the principal
eigenvalues  of time independent or periodic random dispersal eigenvalue problems
\eqref{DC-random-eigen}, \eqref{NC-random-eigen} and \eqref{periodic-random-eigen} (see \cite{Cantrell2003}, \cite{Daners1992}, \cite{Hess1991}, and references therein) and  their
nonlocal dispersal counterparts \eqref{DC-eigen}, \eqref{NC-eigen} and \eqref{periodic-eigen} (see
\cite{Coville2010}, \cite{Coville2013A}, \cite{Garcia2009},
\cite{Hetzer2013},
\cite{Kao2012}, \cite{Rawal2012},
\cite{ShenJDE2010}, \cite{ShenDCDS2015}, \cite{ShenJDE2015} and
references therein).
It is known that a random dispersal operator has always a principal eigenvalue, but a nonlocal dispersal operator may not have a
principal eigenvalue (see \cite{Coville2010} and \cite{ShenJDE2010} for some examples), which reveals some essential
differences between nonlocal dispersal and random dispersal
operators.
 We point out that, very recently, Ding and
Liang \cite{Ding2015}  studied the existence of the
principal eigenvalues of generalized convolution or integral operators on the
circle in periodic media. They showed that such  a generalized convolution
operator has a principal eigenvalue provided that its kernel satisfies the so called uniformly irreducible condition
( see \cite[Theorem
2.4]{Ding2015}).

When $K>1$, principal eigenvalue theory for  \eqref{DC-random-eigen}, \eqref{NC-random-eigen} and
\eqref{periodic-random-eigen} has also been well established. For example, it is known that the principal eigenvalue of \eqref{DC-random-eigen}, (resp. \eqref{NC-random-eigen},
\eqref{periodic-random-eigen}) always exists.
However,
to our best knowledge, the principal eigenvalue problems \eqref{DC-eigen},
\eqref{NC-eigen} and \eqref{periodic-eigen} associated
to coupled nonlocal dispersal systems are hardly studied in
the literature. In this paper, we will  establish some criteria for
the existence of principal eigenvalues of \eqref{DC-eigen},
\eqref{NC-eigen} and \eqref{periodic-eigen}.  The results obtained in this paper
extend many existing ones for
principal eigenvalues of time independent and time periodic
 nonlocal dispersal equations  to time periodic
coupled nonlocal dispersal systems. They provide basic tools for the study of
nonlinear nonlocal dispersal evolution systems.
 We will discuss some applications of the results established in
this paper elsewhere.

The rest of the paper is organized as follows. In Section 2 we
introduce standing notations and state the main results of this
paper. In Section 3, we present some basic properties for coupled nonlocal dispersal systems. In Section 4, we prove the main results.

\vspace{-0.1in}
\section{Notations and Main Results}

\vspace{-0.1in}
In this section, we first introduce the standing notations, present a comparison principle, and give the definition
of principal spectrum points and principal eigenvalues of
\eqref{DC-eq}, \eqref{NC-eq} and \eqref{periodic-eq}, which will be used throughout this paper. We then state the main results of the current paper.

For any $\alpha=(\alpha_{1},...,\alpha_{K})^\top$, we put $|\alpha|=\sqrt{\sum^{K}_{i=1}\alpha^{2}_{i}}$. Let
$$
(\RR^K)^+=\{\alpha=(\alpha_{1},...,\alpha_{K})^\top\,|\, \alpha_i\in\RR,\alpha_i\ge 0,i=1,2,\cdots,K\}
$$
and
$$
(\RR^K)^{++}=\{\alpha=(\alpha_{1},...,\alpha_{K})^\top\,|\, \alpha_i\in\RR,\alpha_i> 0,i=1,2,\cdots,K\}.
$$
Let
\vspace{-0.05in}$$
X_1=X_{2}=C(\bar D,\RR^K)
\vspace{-0.05in}$$
with norm $\|\mathbf{u}\|_{X_{i}}=\sup_{x\in \bar D}|\mathbf{u}(x)|\
(i=1,2)$,
$$
X_3=\{\mathbf{u}\in C_{\rm unif}^b(\RR^N,\RR^K)|
\mathbf{u}(x+p_l\mathbf{e}_l)=\mathbf{u}(x),\ l=1,...,N\}
$$
with norm $\|\mathbf{u}\|_{X_{3}}=\sup_{x\in
\Bbb{R}^{N}}|\mathbf{u}(x)|$, and
$$
X_i^+=\{\mathbf{u}\in X_i\,|\, \mathbf{u}(x)\in(\RR^K)^+\},\quad X_i^{++}=\{\mathbf{u}\in X_i\,|\, \mathbf{u}(x)\in(\RR^K)^{++}\},\,\,\, i=1,2,3.
\vspace{-0.05in}$$
For given $\mathbf{u},\mathbf{v}\in X_i$, we define
\vspace{-0.1in}$$
\mathbf{u}\le \mathbf{v}\,\,\, (\mathbf{v}\ge \mathbf{u})\quad {\rm if}\quad \mathbf{v}-\mathbf{u}\in X_i^+
\quad {\rm and}\quad
\mathbf{u}\ll \mathbf{v}\,\,\, (\mathbf{v}\gg \mathbf{u})\quad {\rm if}\quad \mathbf{v}-\mathbf{u}\in X_i^{++}.
\vspace{-0.1in}$$

Let
$$
\mathcal{X}_1=\mathcal{X}_2=\{\mathbf{u}\in C(\RR\times\bar
D,\RR^K)\,|\, \mathbf{u}(t+T,x)=\mathbf{u}(t,x)\}
$$
with norm $\|\mathbf{u}\|_{\mathcal{X}_{i}}=\sup_{t\in\Bbb{R},x\in \bar
D}|\mathbf{u}(t,x)|\ (i=1,2)$,
$$
\mathcal{X}_3=\{\mathbf{u}\in C(\RR\times\RR^N,\RR^K)\,|\,\mathbf{
u}(t+T,x)=\mathbf{u}(t,x+p_l\mathbf{e}_l)=\mathbf{u}(t,x), \
l=1,...,N\}
$$
with norm $\|\mathbf{u}\|_{\mathcal{X}_{3}}=\sup_{t\in\Bbb{R},x\in
\Bbb{R}^{N}}|\mathbf{u}(t,x)|$, and
\vspace{-0.05in}\[
\mathcal{X}^{+}_{i}=\{\mathbf{u}\in \mathcal{X}_{i}|\ \mathbf{u}(t,x)\in (\RR^K)^{+}\},\quad \mathcal{X}^{++}_{i}=\{\mathbf{u}\in \mathcal{X}_{i}|\ \mathbf{u}(t,x)\in (\RR^K)^{++}\},\,\, i=1,2,3.
\vspace{-0.1in}\]
For any $\mathbf{u},\mathbf{v}\in \mathcal{X}_{i}$, we define
\vspace{-0.05in}\[
\mathbf{u}\leq \mathbf{v}\ (\mathbf{v}\geq
\mathbf{u})\quad\text{if}\quad \mathbf{v}-\mathbf{u}\in
\mathcal{X}^{+}_{i}\,\,\,{\rm and}\,\,\, \mathbf{u}\ll \mathbf{v}\ (\mathbf{v}\gg
\mathbf{u})\quad\text{if}\quad \mathbf{v}-\mathbf{u}\in
\mathcal{X}^{++}_{i}.
\vspace{-0.05in}\]

General semigroup theory (see \cite{Pazy1983}) guarantees that for every
$\mathbf{u}_{0}\in X_{i}$, \eqref{DC-eq} (resp, \eqref{NC-eq},
\eqref{periodic-eq}) has a unique solution
$\mathbf{u}_{1}(t,x;s,\mathbf{u}_{0})$ (resp,
$\mathbf{u}_{2}(t,x;s,\mathbf{u}_{0})$,
$\mathbf{u}_{3}(t,x;s,\mathbf{u}_{0})$) with
$\mathbf{u}_{1}(s,\cdot;s,\mathbf{u}_{0})=\mathbf{u}_{0}\in X_{1}$
(resp, $\mathbf{u}_{2}(s,\cdot;s,\mathbf{u}_{0})=\mathbf{u}_{0}\in
X_{2}$, $\mathbf{u}_{3}(s,\cdot;s,\mathbf{u}_{0})=\mathbf{u}_{0}\in
X_{3}$). Put
\[
(\mathbf{\Phi}_{i}(t,s;A_{i})\mathbf{u}_{0})(\cdot):=\mathbf{u}_{i}(t,\cdot;s,\mathbf{u}_{0}),\quad
i=1,2,3.
\]

\vspace{-0.15in}
\begin{definition}\label{def2.1}
A continuous function $\mathbf{u}(t,x)$ on $[0,\tau)\times \bar D$
is called a supersolution (or subsolution) of \eqref{DC-eq} if for
any $x\in \bar D$, $\mathbf{u}(t,x)$ is continuously differentiable on $[0,\tau)$
and satisfies that for each $x\in \bar D$,
\vspace{-0.1in}\[
\mathbf{u}_t\geq\, (\text{or }\leq)\int_{D}\kappa
(y-x)\mathbf{u}(t,y)dy-\mathbf{u}(t,x)+A_{1}(t,x)\mathbf{u}(t,x),\quad
x\in\bar D
\vspace{-0.1in}\]
for $t\in[0,\tau)$. Supersolutions and subsolutions of \eqref{NC-eq}
and \eqref{periodic-eq} are defined in an analogous way.
\end{definition}

\vspace{-0.1in}
\begin{proposition}\label{pro-comparison}{\rm (Comparison
Principle)}
Assume that (H1) and (H2) hold.
\begin{itemize}
\vspace{-0.05in}\item[(1)] If $\mathbf{u}^{-}(t,x)$ and $\mathbf{u}^{+}(t,x)$ are bounded
subsolution and supersolution of \eqref{DC-eq} {\rm (}resp. \eqref{NC-eq},
\eqref{periodic-eq}{\rm )} on $[0,\tau)$, respectively, and
$\mathbf{u}^{-}(0,x)\leq \mathbf{u}^{+}(0,x)$ for any
$x\in D_1$ (resp. $x\in D_2$, $x\in D_3$), then $ \mathbf{u}^{-}(t,x)\leq
\mathbf{u}^{+}(t,x)$ for any $t\in[0,\tau)$ and $x\in D_1$ (resp. $x\in D_2$, $x\in D_3$).

\vspace{-0.05in}\item[(2)] For given $1\le i\le 3$ and  $\mathbf{u}_{0}\in X^{+}_{i}$,
$\mathbf{u}_{i}(t,\cdot;s,\mathbf{u}_{0})$ exists for all $t\geq s$.

\vspace{-0.05in}\item[(3)]  For given $1\le i\le 3$ and $\mathbf{u}^{1}_{0}(\cdot), \mathbf{u}^{2}_{0}(\cdot)\in X_{i}$,
if $\mathbf{u}^{1}_{0}\leq \mathbf{u}^{2}_{0}$ and
$\mathbf{u}^{1}_{0}\not\neq \mathbf{u}^{2}_{0}$, then
$\mathbf{\Phi}_{i}(t,s;A_{i})\mathbf{u}^{1}_{0}\ll\mathbf{\Phi}_{i}(t,s;A_{i})\mathbf{u}^{2}_{0}$
for all $t>s$.
\end{itemize}
\end{proposition}

\begin{proof}
Note that $A_{i}(t,x)$ is  cooperative and strongly
irreducible for any $t\in\Bbb{R}$ and $x\in D_{i}$. Then the results
of (1)-(3) follows from the similar argument in \cite[Proposition
2.1 and 2.2]{ShenJDE2010} and \cite[Proposition 3.1 and
3.2]{Hetzer2012}.
\end{proof}

\vspace{-0.1in} Define
\vspace{-0.05in}\begin{align*}
(\mathcal{K}_i\mathbf{u})(t,x):=&\int_{D}\kappa(y-x)\mathbf{u}(t,y)dy
\quad {\rm for}\quad \mathbf{u}\in \mathcal{X}_i,\quad i=1,2\\
(\mathcal{K}_3\mathbf{u})(t,x):=&\int_{\RR^N}\kappa(y-x)\mathbf{u}(t,y)dy\quad{\rm
for}\quad \mathbf{u}\in\mathcal{X}_3
\vspace{-0.05in}\end{align*}
and
\vspace{-0.05in}\begin{align*}
(\mathcal{A}_i
\mathbf{u})(t,x):=&-\mathbf{u}(t,x)+A_{i}(t,x)\mathbf{u}(t,x)\quad {\rm for}\,\, \mathbf{u}\in \mathcal{X}_i,\,\, i=1,3\\
(\mathcal{A}_2 \mathbf{u})(t,x):=&- \int_D\kappa(y-x)dy \cdot
\mathbf{u}(t,x)+A_{2}(t,x)\mathbf{u}(t,x) \quad {\rm for}\,\, \mathbf{u}\in \mathcal{X}_2.
\end{align*}
\begin{definition}\label{def2.2}
For given $1\le i\le 3$, let $\sigma(-\p_t +\mathcal{K}_i+\mathcal{A}_i)$ be the spectrum of
$-\p_t +\mathcal{K}_i+\mathcal{A}_i$ on $\mathcal{X}_{i}$. Let
\[
\lambda_{i}(A_{i})=\sup\{{\rm Re}\lambda|\, \lambda\in \sigma(-\p_t
+\mathcal{K}_i+\mathcal{A}_i)\}.
\]
We call $\lambda_{i}(A_{i})$  the {\rm principal spectrum point} of
$-\p_t +\mathcal{K}_i+\mathcal{A}_i$. If
$\lambda_{i}(A_{i})$ is an isolated eigenvalue of $-\p_t
+\mathcal{K}_i+\mathcal{A}_i$ with a positive eigenfunction
$\mathbf{v}$ (i.e, $\mathbf{v}\in \mathcal{X}^{+}_{i}$), then
$\lambda_{i}(A_{i})$ is called the {\rm principal eigenvalue} of $-\p_t
+\mathcal{K}_i+\mathcal{A}_i$ or it is said that {\rm $-\p_t
+\mathcal{K}_i+\mathcal{A}_i$ has a principal eigenvalue}.
\end{definition}

\vspace{-0.1in}
\begin{lemma}\label{lem2.1}
For given $1\le i\le 3$, assume that $A_{i}(t,x)$ satisfies (H1) and (H2). For any given
$x\in D_{i}$, the eigenvalue problem
\begin{equation}\label{2.1}
\begin{cases}
-\frac{d\boldsymbol{\phi}(t)}{dt}+A_{i}(t,x)\boldsymbol{\phi}(t)=\lambda \boldsymbol{\phi}(t),\quad \boldsymbol{\phi}(t)\in\RR^K\\
\boldsymbol{\phi}(t+T)=\boldsymbol{\phi}(t)
\end{cases}
\end{equation}
has a unique
real eigenvalue, denoted by  $\lambda_{i}(x)$,  which has a unique corresponding positive eigenfunction, denoted by
$\boldsymbol{\phi}_{i}(t,x)$, with $|\boldsymbol{\phi}_i(0,x)|=1$.
\end{lemma}
\begin{proof}
$A_{i}(t,x)$ is  a cooperative,
strongly irreducible, and T-periodic matrix for any given $x\in D_{i}$
by (H1) and (H2).
Let $\mathbf{I}$ be the $K\times K$ identity matrix and
$\mathbf{U}(t;x)$ be the fundamental matrix solution of
$\frac{d}{dt}\mathbf{U}(t;x)=A_{i}(t,x)\mathbf{U}(t;x)$ with
$\mathbf{U}(0;x)=\mathbf{I}$ for any given $x\in D_{i}$. By
\cite[Theorem 4.1.1]{Smith1995},
$\mathbf{U}(T;x):\RR^K\rightarrow \RR^K$ is a compact
and strongly positive linear operator. Then  the Krein-Rutmann
theorem (or the Perron-Frobenius theorem in our present
finite-dimensional case, see \cite[Theorem 4.3.1]{Smith1995}), implies  for
any given $x\in D_{i}$ that
the spectral radius $r(\mathbf{U}(T;x))$ is an algebraic simple isolated eigenvalue of $\mathbf{U}(T;x)$ with
a unique positive eigenvector $\boldsymbol{\phi}^i(x)\in(\RR^K)^+$ with $|\boldsymbol{\phi}^i(x)|=1$.  The lemma follows with
$\lambda_i(x)=\frac{1}{T}\ln r(\mathbf{U}(T;x))$ and  $\boldsymbol{\phi}_i(t,x)=e^{-\lambda_i(x)t}\mathbf{U}(t;x)\boldsymbol{\phi}^i(x)$.
\end{proof}

We remark that $\lambda_i(x)$ and $\boldsymbol{\phi}_i(t,x)$ are as smooth in $x$ as $A_i(t,x)$ in $x$. We also remark that, when
$A_i(t,x)\equiv A_i(x)$, $\lambda_i(x)$ is the largest real part of the eigenvalues of the matrix $A_i(x)$.
Let
\begin{equation}\label{2.2}
h_{i}(x)=
\begin{cases}
-1+\lambda_{i}(x),\quad &\text{for }i=1,3,\\
-\int_{D}k(y-x)dy+\lambda_{2}(x),\quad &\text{for }i=2.
\end{cases}
\end{equation}

Our main results on  the principal spectral points and principal
eigenvalues of coupled nonlocal dispersal operators
 can
then be stated as follows.

\vspace{-0.1in} \begin{theorem}\label{th2.1}
Let $1\leq i\leq 3$ be given.
 $\lambda_{i}(A_{i})$ is the principal eigenvalue of
$-\partial_{t}+\mathcal{K}_{i}+\mathcal{A}_{i}$ if and only if
$\lambda_{i}(A_{i})>\max_{x\in D_{i}}h_{i}(x)$.
\end{theorem}

\begin{theorem}\label{th2.2}
Let $1\leq i\leq 3$ be given.
\begin{itemize}
\vspace{-0.05in}\item[(1)] Let $\alpha_i=\max_{x\in D_i}h_i(x)$. If there is a bounded domain $D_0\subset D_i$ such that
$1/(\alpha_i-h_i(\cdot))\not \in L^1(D_0)$, then the principal eigenvalue of
$-\partial_{t}+\mathcal{K}_{i}+\mathcal{A}_{i}$ exists.

\vspace{-0.05in}\item[(2)] Let $\eta_{i}=\min_{t\in\Bbb{R},x\in D_{i}}\{\phi_{i1}(t,x),...,\phi_{iK}(t,x)\}$
and $\tilde\eta_i =\max_{t\in\Bbb{R},x\in D_{i}}\{\phi_{i1}(t,x),...,\phi_{iK}(t,x)\}$
for $i=1,3$, where $\boldsymbol{\phi}_i(t,x)= ( \phi_{i1}(t,x),...,\phi_{iK}(t,x))^\top$ is as in Lemma \ref{lem2.1}. If $\max_{x\in D_{i}}\lambda_{i}(x)-\min_{x\in
D_{i}}\lambda_{i}(x)<\frac{\eta_{i}}{\tilde \eta_i}\inf_{x\in D_{i}} \int_{D_{i}} \kappa(y-x)dy$ in the case
$i=1$ and $\max_{x\in D_{i}}\lambda_{i}(x)-\min_{x\in
D_{i}}\lambda_{i}(x)<\frac{\eta_{i}}{\tilde\eta_i}$ in the case $i=3$, then the principal
eigenvalue of $-\partial_{t}+\mathcal{K}_{i}+\mathcal{A}_{i}$
exists for $i=1,3$.

\vspace{-0.05in}\item[(3)]  Suppose that $\kappa(z)=\frac{1}{\delta^{N}}\widetilde{\kappa}(\frac{z}{\delta})$  for some $\delta>0$ and $\widetilde{\kappa}(\cdot)$ with $\widetilde{\kappa}(z)\geq 0$, ${\rm supp}(\widetilde{\kappa})=B(0,1):=\{z\in\Bbb{R}^{N}|\ \|z\|<1\}$, $\int_{\Bbb{R}^{N}}\widetilde{\kappa}(z)dz=1$, and $\widetilde{k}(\cdot)$ being symmetric with respect to $0$. Then the principal eigenvalue of $-\partial_{t}+\mathcal{K}_{i}+\mathcal{A}_{i}$
exists for $0<\delta\ll 1$.
\end{itemize}
\end{theorem}

\begin{theorem}\label{th2.3}
Let $1\leq i\leq 3$ be given. If $\lambda_{i}(A_{i})$ is the principal eigenvalue of $-\partial_{t}+\mathcal{K}_{i}+\mathcal{A}_{i}$, then it is an algebraically simple eigenvalue of $-\partial_{t}+\mathcal{K}_{i}+\mathcal{A}_{i}$ with a positive eigenfunction.
\end{theorem}

\begin{corollary}
\label{cor2.3}
Let $1\leq i\leq 3$ be given.  If $h_{i}(x)$ is $C^{N}$ and there is some $x_{0}\in
D_{i}$ such that $h_{i}(x_{0})=\max_{x\in
D_{i}}h_{i}(x)$ and the partial derivatives  of $h_i(x)$ up
to order $N-1$ at $x_{0}$ are zero, then the principal eigenvalue of
$-\partial_{t}+\mathcal{K}_{i}+\mathcal{A}_{i}$ exists.
\end{corollary}

\begin{proof}[Proof of Corollary \ref{cor2.3}]
Assume that  $h_{i}(x)$ is $C^{N}$ and there is some $x_{0}\in
D_{i}$ such that $\alpha_i=h_{i}(x_{0})=\max_{x\in
D_{i}}h_{i}(x)$, and the  partial derivatives of $h_i(x)$  up
to order $N-1$ at $x_{0}$ are zero. Then there are bounded domain $D_0\subset D_i$ and
$M>0$ such that
\vspace{-0.05in}\[
0\le \alpha_i-h_{i}(x)\leq
M |x-x_{0}|^{N},\quad x\in D_0.
\vspace{-0.05in}\]
Therefore
$$
\frac{1}{\alpha_i-h_i(x)}\ge \frac{1}{M|x-x_0|^{N}}\quad {\rm for}\quad x\in D_0.
$$
This implies that $1/(\alpha_0-h_i(\cdot))\not \in L^1(D_0)$. The
corollary then follows from Theorem \ref{th2.2}(1).
\end{proof}

\begin{remark}\label{rem2.1}
{\rm
\begin{itemize}
\item[(1)] The condition in Corollary \ref{cor2.3} is called the {\it vanishing condition}.

\vspace{-0.05in}\item[(2)] It remains open whether the analogous result holds for the case $i=2$ in  Theorem \ref{th2.2}(2).
When $K=1$,  Theorem \ref{th2.2} (2) holds for $i=1,2,3$ under some weaker conditions (see \cite[Theorem B (2)]{Rawal2012}).

\vspace{-0.05in}\item[(3)]  If $N=1$ or $2$, $h_i(\cdot)$ is $C^N$, and $\max_{x\in D_i}h_i(x)=h_i(x_0)$ for some $x_i\in{\rm Int}(D_i)$, then the conditions in Theorem \ref{th2.2} (1) are always satisfied
and hence the principal eigenvalue of $-\p_t+\mathcal{K}_i+\mathcal{A}_i$  exists. The principal eigenvalue may not exist if $N\geq 3$ (see \cite{Rawal2012}, \cite{ShenJDE2010}).
\end{itemize}}
\end{remark}

\vspace{-0.1in}
\section{Basic Properties}

\vspace{-0.1in} In this section, we present some basic properties for principal
spectrum points of coupled nonlocal dispersal systems
\eqref{DC-eigen}, \eqref{NC-eigen} and \eqref{periodic-eigen}.
Throughout this section, we assume that  (H1) and (H2) hold.

Let $\mathcal{K}_{i}:\mathcal{X}_{i}\rightarrow
\mathcal{X}_{i}$ and $\mathcal{A}_{i}:\mathcal{X}_{i}\rightarrow
\mathcal{X}_{i}$ be defined as in Section 2 for $i=1,2,3$.
Consider the eigenvalue problem,
\begin{equation}\label{3.1}
-\partial_{t}\mathbf{u}+\mathcal{K}_{i}\mathbf{u}+\mathcal{A}_{i}\mathbf{u}=\lambda\mathbf{u},\quad
i=1,2,3.
\end{equation}
 Recall that $\mathbf{I}$ denotes the $K\times K$ identity matrix. Without loss of generality, we assume that $-A_{i}(t,x)$ is positive definite for any $(t,x)\in\Bbb{R}\times D_{i}$, that is, for any $\mathbf{u}\in\RR^K$,
\vspace{-0.05in}\begin{equation}\label{3.01}
-\mathbf{u}^{T}A_{i}(t,x)\mathbf{u}\geq 0\quad \forall (t,x)\in \Bbb{R}\times D.
\vspace{-0.05in}\end{equation}
In fact,  there exists $M>0$ such that for any $\mathbf{u}\in\RR^K$,
\vspace{-0.05in}\[
-\mathbf{u}^{T}(A_{i}(t,x)-M\mathbf{I})\mathbf{u}\geq 0\quad \forall (t,x)\in \Bbb{R}\times D.
\vspace{-0.05in}\]
Hence $-\widetilde A_{i}(t,x)$ is positive definite for any $(t,x)\in\RR\times D_i$, where
 $\widetilde A_{i}(t,x)={A}_{i}(t,x)-M\mathbf{I}$. If $-A_{i}(t,x)$ is not positive definite for some $(t,x)\in\RR\times D_i$,
we may consider the eigenvalue problem,
\vspace{-0.05in}\[
-\partial_{t}\mathbf{u}+\mathcal{K}_{i}\mathbf{u}+\widetilde{\mathcal{A}}_{i}\mathbf{u}=\widetilde{\lambda}\mathbf{u},
\vspace{-0.05in}\]
where $\widetilde{\mathcal{A}}_{i}\mathbf{u}={\mathcal{A}}_{i}\mathbf{u}-M\mathbf{u}$.

Let $\mathcal{X}_{i}\oplus
i\mathcal{X}_{i}=\{\mathbf{u}+i\mathbf{v}|\mathbf{u},\mathbf{v}\in
\mathcal{X}_{i}\}$. Observe that  if $\alpha\in
\mathbb{C}$ is such that
$\alpha\mathbf{I}+\partial_{t}-\mathcal{A}_{i}$ is invertible, then
\eqref{3.1} with $\lambda=\alpha$  has nontrivial solutions in $\mathcal{X}_{i}\oplus
i\mathcal{X}_{i}$ is equivalent to
\vspace{-0.05in}\[
\mathcal{K}_{i}(\alpha\mathbf{I}+\partial_{t}-\mathcal{A}_{i})^{-1}\mathbf{v}=\mathbf{v}
\vspace{-0.05in}\]
has nontrivial solutions in $\mathcal{X}_{i}\oplus
i\mathcal{X}_{i}$. Recall that $h_{i}(x)$ is given by \eqref{2.2}.

\begin{proposition}\label{pro3.1}
Let $1\leq i\leq 3$ be given. Then
\[
\big[\min_{x\in D_{i}}h_{i}(x), \max_{x\in
D_{i}}h_{i}(x) \big]\subset
\sigma(-\partial_{t}+\mathcal{A}_{i}).
\]
\end{proposition}
\begin{proof}
Assume that  $h_{i}(x_0)\in \rho(-\partial_{t}+\mathcal{A}_{i})$ for some $x_0\in D_i$. Put $h_i=h_i(x_0)$.
Then
\vspace{-0.05in}\[
-\partial_{t}\mathbf{u}+\mathcal{A}_{i}\mathbf{u}-h_{i}\mathbf{u}=\mathbf{v}(t)
\vspace{-0.05in}\]
has a unique solution $\mathbf{u}(\cdot,\cdot; x_{0},\mathbf{v})\in
\mathcal{X}_{i}$ for any $\mathbf{v}\in \mathcal{X}_{i}$ with
$\mathbf{v}(t,x)\equiv \mathbf{v}(t)$.  This implies that $\mathbf{u}(t;x_0,\mathbf{v}):=\mathbf{u}(t,x_0;x_0,\mathbf{v})$ is a solution of
\vspace{-0.05in}$$
-\p_t \mathbf{u}+A_i(t,x_0) \mathbf{u}-\lambda_i(x_0)\mathbf{u}=\mathbf{v}(t).
\vspace{-0.05in}$$
Therefore the Fredholm alternative implies that
\vspace{-0.05in}\[
-\frac{d\mathbf{u}}{dt}+A_{i}(t,x_{0})\mathbf{u}-\lambda_{i}(x_{0})\mathbf{u}=0
\vspace{-0.05in}\]
has no nontrivial solution $\mathbf{u}(t)$ with $\mathbf{u}(t+T)=\mathbf{u}(t)$.
Recall that $\boldsymbol{\phi}_{i}(t,x_{0})$ is
an eigenfunction of \eqref{2.1} with $x=x_0$  corresponding to
$\lambda_{i}(x_{0})$. Then
$\mathbf{u}(t;x_0,\boldsymbol{\phi}_{i}):=\boldsymbol{\phi}_{i}(t,x_{0})$ satisfies
\[
-\partial_{t}\mathbf{u}(t;x_{0},\boldsymbol{\phi}_{i})+A_{i}(t,x_{0})\mathbf{u}(t;x_{0},\boldsymbol{\phi}_{i})
-\lambda_{i}(x_{0})\mathbf{u}(t;x_{0},\boldsymbol{\phi}_{i})=0,
\]
 which is a contradiction. Hence $h_{i}(x_0)\in
\sigma(-\partial_{t}+\mathcal{A}_{i})$  for any $x_0\in D_i$ and the proposition follows.
\end{proof}

\vspace{-0.05in}
Let
$$
X_0=\begin{cases}C(D_i,\RR),\quad i=1,2\cr
\{u\in C(D_i,\RR)\,|\, u(\cdot+p_l\mathbf{e}_l)=u(\cdot)\},\quad i=3.
\end{cases}
$$

\begin{proposition}\label{pro3.2}
Let $1\leq i\leq 3$ be given, then
$(\alpha\mathbf{I}+\partial_{t}-\mathcal{A}_{i})^{-1}$ exists
for every $\alpha\in \mathbb{C}$ with
${\rm Re\, } \alpha>\max_{x\in D_{i}}h_{i}(x)$.
Moreover, one has for any $\alpha\in\mathbf{C}$ with ${\rm Re}\alpha>\max_{x\in
D_{i}}h_{i}(x)$ and any $v_i(\cdot)\in X_0$ that
\[
\left((\alpha\mathbf{I}+\partial_{t}-\mathcal{A}_{i})^{-1}\boldsymbol{\phi}_{i} v_i\right)(t,x)=
\frac{1}{\alpha-h_{i}(x)}\boldsymbol{\phi}_{i}(t,x)v_i(x)\quad \forall\,\, (t,x)\in \RR\times D_i,
\]
where
$\boldsymbol{\phi}_{i}(t,x)\in \mathcal{X}^{+}_{i}$ is as in Lemma \ref{lem2.1}.
\end{proposition}
\begin{proof}
Let $\widetilde{X}=\{\mathbf{u}(t)\in C(\Bbb{R},\Bbb{R}^{K})|\ \mathbf{u}(t+T)=\mathbf{u}(t)\}$. First of all,   Floquent theory for periodic ordinary differential equations implies that for any $\alpha\in
\mathbb{C}$ with ${\rm Re\,}\alpha>\max_{x\in
D_{i}}h_{i}(x)$ and any given $x\in D_{i}$,
$(\alpha\mathbf{I}+\partial_{t}-\mathcal{A}_{i}(x))^{-1}$ exists in $\widetilde{X}$, where $(\mathcal{A}_i(x)\mathbf{u})(t)=(\mathcal{A}\mathbf{u})(t,x)$
for $\mathbf{u}\in \widetilde{X}$. Hence $(\alpha\mathbf{I}+\partial_{t}-\mathcal{A}_{i})^{-1}$ exists in $\mathcal{X}_i$.

 By Lemma \ref{lem2.1}, $\lambda_{i}(x)$ is the principal eigenvalue of \eqref{2.1} with the positive eigenvector $\boldsymbol{\phi}_{i}(t,x)$  for any $1\leq i\leq 3$ and any given $x\in D_{i}$. By the definition of $\mathcal{A}_{i}$ and $h_{i}(x)$,  we have
\[
-\partial_{t}\boldsymbol{\phi}_{i}(t,x)+(\mathcal{A}_{i}\boldsymbol{\phi}_{i})(t,x)=h_{i}(x)\boldsymbol{\phi}_{i}(t,x).
\]
for any $x\in D_{i}$, which
implies that
\[
-\partial_{t}\boldsymbol{\phi}_{i}(t,x)v_i(x)+(\mathcal{A}_{i}\boldsymbol{\phi}_{i} v_i)(t,x)=h_{i}(x)\boldsymbol{\phi}_{i}(t,x)v_i(x).
\]
for every $v_i\in X_0$ and given $x\in D_i$, hence
\[
\Big((\alpha\mathbf{I}+\partial_{t}-\mathcal{A}_{i})\boldsymbol{\phi}_{i}v_i\Big)(t,x)
=(\alpha-h_{i}(x))\boldsymbol{\phi}_{i}(t,x)v_i(x)
\]
for every $\alpha\in\mathbf{C}$, each $x\in D_{i}$ and each $v_i\in X_0$.
It then follows that for any $\alpha\in \mathbf{C}$ with ${\rm Re\, } \alpha>\max_{x\in D_{i}}h_{i}(x)$ and any $v_i\in X_0$,
\vspace{-0.1in}\begin{equation}\label{3.002}
\Big((\alpha\mathbf{I}+\partial_{t}-\mathcal{A}_{i})^{-1}\boldsymbol{\phi}_{i} v_i\Big)(t,x)
=\frac{1}{\alpha-h_{i}(x)}\boldsymbol{\phi}_{i}(t,x) v_i(x)\quad \forall \,\, (t,x)\in\RR\times D_i.
\vspace{-0.1in} \end{equation}
   This completes the proof.
\end{proof}

\vspace{-0.1in} \begin{proposition}\label{pro3.3}
Let $1\leq i\leq 3$ be given.
\begin{itemize}
\vspace{-0.1in} \item[(1)] $-\partial_{t}+\mathcal{A}_{i}$ generates a
positive semigroup of contractions on $\mathcal{X}_{i}$.

 \vspace{-0.1in}\item[(2)]
$\mathcal{K}_{i}(\alpha\mathbf{I}+\partial_{t}-\mathcal{A}_{i})^{-1}$
is  compact for every
$\alpha\in \mathbb{C}$ with ${\rm
Re}\alpha>\max_{x\in D_{i}}h_{i}(x)$.
\end{itemize}
\end{proposition}
\begin{proof}
(1) Define $T(s)\mathbf{u}$ by
\vspace{-0.1in}\[
(T(s)\mathbf{u})(t,x)=\mathbf{u}(t-s,x).
\vspace{-0.1in}\]
for $s\ge 0$ and $\mathbf{u}\in \mathcal{X}_{i}$.
Then $T(s)$ is a continuous semigroup generated by  $-\partial_{t}$ on $\mathcal{X}_{i}$. It is clear that $\|T(s)\mathbf{u}\|\le \|\mathbf{u}\|$ for all
$s\ge 0$ and $\mathbf{u}\in \mathcal{X}_i$, and that for any $\mathbf{u}\in \mathcal{X}_{i}$ with $\mathbf{u}(\cdot,\cdot)\geq \mathbf{0}$, $T(s)\mathbf{u}\geq \mathbf{0}$ for any $s\geq 0$.  Therefore, $\{T(s)\}_{s\in\Bbb{R}^{+}}$ is a positive continuous semigroup of contraction on $\mathcal{X}_{i}$ with generator  $-\partial_{t}$.

On the other hand, given  $\mu>0$ and $\mathbf{u}\in \mathcal{X}_{i}$,  \eqref{3.01} yields
\[
\|\mathbf{u}\| \|(\mu\mathbf{I}-\mathcal{A}_{i})\mathbf{u}\|\geq |\mathbf{u}^{T}(t,x)((\mu\mathbf{I}-\mathcal{A}_{i})\mathbf{u})(t,x)|= \mu |\mathbf{u}(t,x)|^2-\mathbf{u}^{T}(t,x)\mathcal{A}_{i}(t,x)\mathbf{u}(t,x)\ge \mu|\mathbf{u}(t,x)|^2
\]
for all $(t,x)\in \Bbb{R}\times D_{i}$,
 which implies that $\|(\mu\mathbf{I}-\mathcal{A}_{i})\mathbf{u}\|\geq \mu\|\mathbf{u}\|$ for any $\mathbf{u}\in \mathcal{X}_{i}$ and $\mu>0$. By \cite[Theorem 1.4.2]{Pazy1983}, we have that $\mathcal{A}_{i}$ is dissipative on $\mathcal{X}_i$. It then follows from \cite[Corollary 3.3.3]{Pazy1983}
that   $-\partial_{t}+\mathcal{A}_{i}$  is the infinitesimal generator of continuous semigroup of contraction on $\mathcal{X}_{i}$.

(2) First, Proposition \ref{pro3.2} guarantees that
$(\alpha\mathbf{I}+\partial_{t}-\mathcal{A}_{i})^{-1}$ exists for every $\alpha$ with ${\rm Re\,}\alpha>\max_{x\in D_{i}}h_{i}(x)$.
Next for any bounded sequence $\{\mathbf{u}_{n}\}\in \mathcal{X}_{i}\oplus\mathcal{X}_{i}$, let
\vspace{-0.05in}\[
\mathbf{w}_{n}=\left(\mathcal{K}_{i}(\alpha\mathbf{I}+\partial_{t}-\mathcal{A}_{i})^{-1}\right)\mathbf{u}_{n}
\vspace{-0.05in}\]
and
\vspace{-0.05in}\[
\mathbf{v}_{n}=\left(\alpha\mathbf{I}+\partial_{t}-\mathcal{A}_{i}\right)^{-1}\mathbf{u}_{n}.
\vspace{-0.05in}\]
Then
\vspace{-0.05in}\[
\mathbf{w}_{n}=\mathcal{K}_{i}\mathbf{v}_{n}.
\vspace{-0.05in}\]
By the boundedness of $\mathcal{A}_{i}$, both $\{\mathbf{v}_{n}\}$ and $\{\partial_{t}\mathbf{v}_{n}\}$ are bounded sequences in $\mathcal{X}_{i}\oplus\mathcal{X}_{i}$.
 Then $\{\mathbf{w}_{n}\}=\{\mathcal{K}_{i}\mathbf{v}_{n}\}$, $\{\partial_{t}\mathbf{w}_{n}\}$, and $\{\partial_{x_i}\mathbf{w}_{n}\}$
 ($i=1,2,\cdots,N$) are bounded sequences in $\mathcal{X}_{i}\oplus\mathcal{X}_{i}$.
 We can then show that $\{\mathbf{w}_{n}\}$ has a convergent subsequence. This completes the proof.
\end{proof}

\begin{proposition}\label{pro3.4}
For given $1\leq i\leq 3$ and $\alpha>\max_{x\in D_i}h_i(x)$, let
$r(\mathcal{K}_{i}(\alpha\mathbf{I}+\partial_{t}-\mathcal{A}_{i})^{-1})$
be the spectral radius of
$\mathcal{K}_{i}(\alpha\mathbf{I}+\partial_{t}-\mathcal{A}_{i})^{-1}$.  Let $x_{0}\in D_{i}$ be such that $h_{i}(x_{0})=\max_{x\in D_i}h_i(x)$. If there is a bounded domain $D_0\subset D_i$ such that
$1/(h_{i}(x_{0})-h_i(\cdot))\not \in L^1(D_0)$,  then there is $\alpha>\max_{x\in D_{i}}h_{i}(x)$
such that
\[
r(\mathcal{K}_{i}(\alpha\mathbf{I}+\partial_{t}-\mathcal{A}_{i})^{-1})>1.
\]
\end{proposition}

\begin{proof}
By Proposition \ref{pro3.2},  $(\alpha\mathbf{I}+\partial_{t}-\mathcal{A}_{i})^{-1}$ exists
for all $\alpha\in\RR$ with $\alpha>h_{i}(x_{0})=\max_{x\in
D_{i}}h_{i}(x)$, and 
\begin{equation}\label{3.2}
\left((\alpha\mathbf{I}+\partial_{t}-\mathcal{A}_{i})^{-1}\boldsymbol{\phi}_{i} v_i\right)(t,x)=
\frac{1}{\alpha-h_{i}(x)}\boldsymbol{\phi}_{i}(t,x) v_i(x)
\end{equation}
for all $v_i\in X_0$.
Let $\delta>0$ be such that $B(0,\delta)\subset {\rm supp}(\kappa(\cdot))$. By \eqref{3.2} we have
\[
\Big(\mathcal{K}_{i}(\alpha\mathbf{I}+\partial_{t}-\mathcal{A}_{i})^{-1}\boldsymbol{\phi}_{i} v_i\Big)(t,x)= \int_{D_{i}}\frac{\kappa(y-x)\boldsymbol{\phi}_{i}(t,y)v_{i}(y)}{\alpha-h_{i}(y)}dy.
\]
Note that   $\boldsymbol{\phi}_{i}(t,x)\geq \eta_{i}\mathbf{e}$ for any $t\in\Bbb{R}$ and $x\in D_{i}$, where
$\eta_i=\min_{t\in\RR,x\in D_i}\{\phi_{i1}(t,x),\cdots,\phi_{iK}(t,x)\}$ and $\mathbf{e}=(1,1,\cdots,1)^\top$.
Then, given $v_i\in X_0$ with $v_i(x)\ge 0$,
\begin{align}\label{3.3}
\Big(\mathcal{K}_{i}(\alpha\mathbf{I}+\partial_{t}-\mathcal{A}_{i})^{-1}\boldsymbol{\phi}_{i} v_i\Big)(t,x)
\geq &
\int_{D_{i}}\frac{\eta_{i} \kappa(y-x) v_i(y)}{\alpha-h_i(y)} dy   \, \mathbf{e}.
\end{align}
Let $\sigma>0$ and $x_1\in D_i$  be such that $2\sigma <\delta$,  $B(x_1,\sigma)\subset D_0$, $B(x_1,2\sigma)\subset D_i$, and
$$
\lim_{\alpha\to h_i(x_0)}\int_{B(\sigma,x_1)}\frac{1}{\alpha-h_i(y)}dy=\infty
$$
($B(x_1,r)=\{x\in\RR^N\,|\, |x|<r\}$).
  Let $v_{i}(x)\in X_{0}$ be such that $v_{i}(x)=1$ if $x\in B(\sigma,x_{1})$ and $v_{i}(x)=0$ if $x\in D_{i}\setminus B(2\sigma,x_{1})$.

Clearly, for every $x\in D_{i}\setminus B(2\sigma,x_{1})$ and $\gamma>1$,
\[
\Big(\mathcal{K}_{i}(\alpha\mathbf{I}+\partial_{t}-\mathcal{A}_{i})^{-1}\boldsymbol{\phi}_{i} v_i\Big)(t,x)>\gamma \boldsymbol{\phi}_{i}(t,x) v_i(x)=0.
\]
For $x\in B(2\sigma,x_{1})$, there is $\widetilde{M}>0$ such that
$\kappa(y-x)\geq \widetilde{M}$ for $y\in B(\sigma,x_{1})$.
It then follows from \eqref{3.3} that for $x\in B(2\sigma,x_{1})$,
\begin{align*}
\Big(\mathcal{K}_{i}(\alpha\mathbf{I}+\partial_{t}-\mathcal{A}_{i})^{-1}\boldsymbol{\phi}_{i} v_i\Big)(t,x)
\geq &
\int_{B(\sigma,x_{1})}\frac{\eta_{i} \kappa(y-x) v_i(y)}{\alpha-h_i(y)} dy   \, \mathbf{e}\\
\geq& \int_{B(\sigma,x_{1})}\frac{1}{\alpha-h_i(y)} dy   \, \eta_{i}\widetilde{M}\mathbf{e}
\end{align*}
Then there is $\gamma>1$ such that, for $\alpha>h_{i}(x_{0})$ and   $\alpha-h_{i}(x_{0})\ll 1$,
\begin{align*}
\Big(\mathcal{K}_{i}(\alpha\mathbf{I}+\partial_{t}-\mathcal{A}_{i})^{-1}\boldsymbol{\phi}_{i} v_i\Big)(t,x)
\geq &\,  \gamma \phi_i(t,x)v_i(x).
\end{align*}
This implies that $r(\mathcal{K}_{i}(\alpha\mathbf{I}+\partial_{t}-\mathcal{A}_{i})^{-1})\ge \gamma>1$
for $\alpha>h_{i}(x_{0})$ and   $\alpha-h_{i}(x_{0})\ll 1$.
\end{proof}



\section{Principal Eigenvalues of Coupled Nonlocal  Dispersal}
\noindent

In this section, we investigate criteria for  the  existence  of principal
eigenvalues of  coupled nonlocal dispersal systems with time
periodic dependence and prove Theorems \ref{th2.1}-\ref{th2.3}. Throughout this section, we assume that (H1) and (H2) hold. Recall that $\lambda_{i}(A_{i})=\sup\{{\rm Re\,}\lambda |\, \lambda\in \sigma(-\partial_{t}+\mathcal{K}_{i}+\mathcal{A}_{i}) \}$ is the principal spectrum point of $-\partial_{t}+\mathcal{K}_{i}+\mathcal{A}_{i}$  ($i=1,2,3$). We first prove a lemma.

\vspace{-0.1in} \begin{lemma}\label{lm4.1}
\begin{itemize}
\item[(1)] For given $1\leq i\leq 3$, if there is $\alpha_{0}>\max_{x\in
D_{i}}h_{i}(x)$ such that
$r(\mathcal{K}_{i}(\alpha_{0}\mathbf{I}+\partial_{t}-\mathcal{A}_{i})^{-1})>1$, then $\lambda_{i}(A_{i})>\max_{x\in
D_{i}}h_{i}(x)$,
$r(\mathcal{K}_{i}(\lambda_{i}(A_{i})\mathbf{I}+\partial_{t}-\mathcal{A}_{i})^{-1})=1$,
and $\lambda_{i}(A_{i})$ is an  isolated eigenvalue
of $-\partial_{t}+\mathcal{K}_{i}+\mathcal{A}_{i}$ of finite multiplicity with a positive
eigenfunction.

\vspace{-0.05in}\item[(2)] If $\lambda_{i}(A_{i})$ is an eigenvalue of $-\partial_{t}+\mathcal{K}_{i}+\mathcal{A}_{i}$ with a positive eigenfunction, then it is geometrically simple.
 \end{itemize}
\end{lemma}

\begin{proof}
(1) Suppose that there is $\alpha_{0}>\max_{x\in D_{i}}h_{i}(x)$ such that $r(\mathcal{K}_{i}(\alpha\mathbf{I}+\partial_{t}-\mathcal{A}_{i})^{-1})>1$.
By Proposition \ref{pro3.3},
$\mathcal{K}_{i}(\alpha\mathbf{I}+\partial_{t}-\mathcal{A}_{i})^{-1}$
is a compact operator for any $\alpha\in\mathbb{C}$
with ${\rm Re\,}\alpha>\max_{x\in
D_{i}}h_{i}(x)$.
It then follows from \cite[Theorem 2.2]{Burger1988} that
$\lambda_{i}(A_{i})>\max_{x\in D_{i}}h_{i}(x)$,
$r(\mathcal{K}_{i}(\lambda_{i}(A_{i})\mathbf{I}+\partial_{t}-\mathcal{A}_{i})^{-1})=1$,
and $\lambda_{i}(A_{i})$ is an  isolated eigenvalue of
$-\partial_{t}+\mathcal{K}_{i}+\mathcal{A}_{i}$ with  finite multiplicity and  a positive
eigenfunction.

(2) We show that $\lambda_{i}(A_{i})$
is geometrically simple. Suppose that $\mathbf{v}(t,x)$ is a
positive eigenfunction of $-\partial_{t}+\mathcal{K}_{i}+\mathcal{A}_{i}$ associated to
$\lambda_{i}(A_{i})$. By Proposition \ref{pro-comparison},
$\mathbf{v}(t,x)>\mathbf{0}$ for $t\in\Bbb{R}$ and $x\in D_{i}$.
Assume that $\mathbf{v}^{\prime}(t,x)$ is also an eigenfunction of
$-\partial_{t}+\mathcal{K}_{i}+\mathcal{A}_{i}$ associated to $\lambda_{i}(A_{i})$. Then there is
$\eta\in\Bbb{R}$  such
that
\[
\mathbf{w}(t,x)\geq 0\quad \forall\,\,  t\in \Bbb{R}, \,\, x\in D_{i}\quad \text{and }\quad \mathbf{w}(t_{0},x_{0})=0,
\]
where $\mathbf{w}(t,x)=\mathbf{v}(t,x)-\eta\mathbf{v}^{\prime}(t,x)$.
Then by Proposition \ref{pro-comparison} again,
$\mathbf{w}(t,x)\equiv 0$ and hence
$\mathbf{v}(t,x)=\eta\mathbf{v}^{\prime}(t,x)$. This implies that
$\lambda_{i}(A_{i})$ is geometrically simple. This completes the
proof.
\end{proof}

\begin{proof}[Proof of Theorem \ref{th2.1}]
We prove the case $i=1$. Other cases can
be proved similarly.

First, assume that $\lambda=\lambda_1(A)$ is  the principal eigenvalue of
$-\partial_{t}+\mathcal{K}_{1}+\mathcal{A}_{1}$ with a positive
eigenfunction $\boldsymbol{\psi}_{1}(t,x)\in\mathcal{X}^{+}_{1}$.
Then  we have
\vspace{-0.05in}\begin{equation}\label{4.1}
-\partial_{t}\boldsymbol{\psi}_{1}(t,x)+(\mathcal{K}_{1}\boldsymbol{\psi}_{1})(t,x)+(\mathcal{A}_{1}\boldsymbol{\psi}_{1})(t,x)
=\lambda_1(A_1)\boldsymbol{\psi}_{1}(t,x),\quad x\in D_1.
\vspace{-0.05in}\end{equation}
Recall that $\mathbf{\Phi}_{1}(t,s;A_{i})$ is the solution operator
of \eqref{DC-eq}. Note
that $\mathbf{u}(t,x)=e^{\lambda t}\boldsymbol{\psi}_{1}(t,x)$ is a
solution of \eqref{DC-eq}, which implies that
\vspace{-0.05in}\[
\left(\mathbf{\Phi}_{1}(t,0;A_{1})\boldsymbol{\psi}_{1}(0,\cdot)\right)(t,x)=e^{\lambda
t}\boldsymbol{\psi}_{1}(t,x).
\vspace{-0.05in}\]
Then by the comparison principle, we have
$\boldsymbol{\psi}_{1}\gg\mathbf{0}$.

 For any given $x_{0}\in D_{1}$, we know that
$\boldsymbol{\phi}_{1}(t,x_{0})$ is an eigenfunction of
\eqref{2.1} corresponding to the eigenvalue $\lambda_{1}(x_{0})$, that is,
\vspace{-0.05in}\begin{equation}\label{4.2}
-\frac{d}{dt}\boldsymbol{\phi}_{1}(t,x_{0})+A_{1}(t,x_{0})\boldsymbol{\phi}_{1}(t,x_{0})=\lambda_{1}(x_{0})\boldsymbol{\phi}_{1}(t,x_{0}).
\vspace{-0.05in}\end{equation}
By \eqref{4.1},
\vspace{-0.05in}\begin{equation}\label{4.3}
-\partial_{t}\boldsymbol{\psi}_{1}(t,x_{0})-\boldsymbol{\psi}_{1}(t,x_{0})+(\mathcal{A}_{1}\boldsymbol{\psi}_{1})(t,x_{0})<\lambda_{1}(A_{1})
\boldsymbol{\psi}_{1}(t,x_{0})
\vspace{-0.05in}\end{equation}
for every $x_{0}\in D_{1}$.  \eqref{4.2} and \eqref{4.3} yield
\vspace{-0.05in}\[
\lambda_{1}(A_{1})>-1+\lambda_{1}(x_{0}) \qquad \mathrm{for }\ x_{0}\in D_{1}.
\vspace{-0.05in}\]
Hence $\lambda_{1}(A_{1})>\max_{x\in D_{1}}\{-1+\lambda_{1}(x)\}$.

Conversely, assume $\lambda_{1}(A_{1})>\max_{x\in
D_{1}}h_{1}(x)$. By Lemma \ref{lm4.1},
$\lambda_{1}(A_{1})$ is the principal eigenvalue of
$-\partial_{t}+\mathcal{K}_{1}+\mathcal{A}_{1}$. This completes the
proof.
\end{proof}

\begin{proof}[Proof of Theorem \ref{th2.2}]
(1)  By Proposition \ref{pro3.4}, $r(\mathcal{K}_{i}(\alpha\mathbf{I}+\partial_{t}-\mathcal{A}_{i})^{-1})>1$
for $0<\alpha-\max_{x\in D_{i}}h_{i}(x)\ll 1$.
Lemma \ref{lm4.1} implies that $\lambda_{i}(A_{i})$ is the principal
eigenvalue of $-\partial_{t}+\mathcal{K}_{i}-\mathcal{A}_{i}$.

(2) We prove the case $i=3$. Case $i=1$ can be proved similarly.  By Proposition \ref{pro3.2}, we have
\[
\left(\mathcal{K}_{3}(\alpha\mathbf{I}+\partial_{t}-\mathcal{A}_{3})^{-1}\boldsymbol{\phi}_{3}v_{3}\right)(t,x)
=\int_{\Bbb{R}^{N}}\frac{\kappa(y-x)}{\alpha-h_{3}(y)}\boldsymbol{\phi}_{3}(t,y)v_{3}(y)dy
\]
for every $\alpha>\max_{x\Bbb{R}^{N}}h_{3}(x)$ and any $v_{3}(x)\in X_{0}$ with $v_{3}(x)\geq 0$.

Note  that $\boldsymbol{\phi}_{3}(t,x)\geq \eta_{3}\mathbf{e}$ for any $t\in\Bbb{R}$ and $x\in\Bbb{R}^{N}$. Let $v_{3}(x)\equiv1$ for all $x\in\Bbb{R}^{N}$. Then
\[
\left(\mathcal{K}_{3}(\alpha\mathbf{I}+\partial_{t}-\mathcal{A}_{3})^{-1}\boldsymbol{\phi}_{3}\right)(t,x)
\geq\int_{\Bbb{R}^{N}}\frac{\kappa(y-x)}{\alpha-h_{3}(y)}dy\ \eta_{3}\mathbf{e}.
\]
Take $\alpha=-1+\max_{x\in \Bbb{R}^{N}}\lambda_{3}(x)+\epsilon$ for $\epsilon\in (0,1)$. Then for any $0<\epsilon<1$,
\begin{align}\label{4.4}
\left(\mathcal{K}_{3}(\alpha\mathbf{I}+\partial_{t}-\mathcal{A}_{3})^{-1}\boldsymbol{\phi}_{3}\right)(t,x)
\geq&\frac{1}{\max_{x\in \Bbb{R}^{N}}\lambda_{3}(x)-\min_{x\in \Bbb{R}^{N}}\lambda_{3}(x)+\epsilon}\int_{\Bbb{R}^{N}}\kappa(y-x)dy\ \eta_{3}\mathbf{e}\nonumber\\
=&\frac{\eta_{3}\mathbf{e}}{\max_{x\in \Bbb{R}^{N}}\lambda_{3}(x)-\min_{x\in \Bbb{R}^{N}}\lambda_{3}(x)+\epsilon}.
\end{align}
  By $\max_{x\in \Bbb{R}^{N}}\lambda_{3}(x)-\min_{x\in \Bbb{R}^{N}}\lambda_{3}(x) < \frac{\eta_{3}}{\tilde\eta_3}$ and  \eqref{4.4},  for  $0<\epsilon\ll 1$, there is $\gamma>1$ such that $\left(\mathcal{K}_{3}(\alpha\mathbf{I}+\partial_{t}-\mathcal{A}_{3})^{-1}\boldsymbol{\phi}_{3}\right)(t,x)
\geq \gamma \boldsymbol{\phi}_{3}(t,x)$ and hence $r(\mathcal{K}_{3}(\alpha\mathbf{I}+\partial_{t}-\mathcal{A}_{3})^{-1})>1$. It then follows from Lemma \ref{lm4.1} that the principal
eigenvalue of $-\partial_{t}+\mathcal{K}_{3}-\mathcal{A}_{3}$ exists.

(3) It can be proved by similar arguments as in \cite[Theorem 2.6]{Kao2012}.
 \end{proof}

\begin{proof}[Proof of Theorem \ref{th2.3}]
Note that for $\alpha>\max_{x\in D_{i}}h_{i}(x)$, $(\alpha\mathbf{I}+\partial_{t}-\mathcal{A}_{i})^{-1}$ exists. For given  $\alpha>\max_{x\in D_{i}}h_{i}(x)$, let
\vspace{-0.05in}\[
(\mathbf{U}_{\alpha}\mathbf{u})(t,x)=\big(\mathcal{K}_{i}(\alpha\mathbf{I}+\partial_{t}-\mathcal{A}_{i})^{-1}\mathbf{u}\big)(t,x)
\quad\text{and}\quad r(\alpha)=r(\mathbf{U}_{\alpha}).
\vspace{-0.05in}\]
By Proposition \ref{pro3.3}, $\mathbf{U}_{\alpha}:\mathcal{X}_{i}\rightarrow \mathcal{X}_{i}$ is a positive and compact operator. Let $\lambda_{i}(A_{i})$ be the principal eigenvalue of $-\partial_{t}+\mathcal{K}_{i}+\mathcal{A}_{i}$. By Lemma \ref{lm4.1}, $\lambda_{i}(A_{i})$ is an isolated geometrically simple eigenvalue of $-\partial_{t}+\mathcal{K}_{i}+\mathcal{A}_{i}$. Let $\alpha_{0}=\lambda_{i}(A_{i})$. Then $r(\alpha_{0})=1$ and $r(\alpha_{0})$ is an isolated geometrically simple eigenvalue of $\mathbf{U}_{\alpha_{0}}$ with $\boldsymbol{\phi}_{i}(\cdot,\cdot)$ being a positive eigenfunction.

The rest of proof is similar to that in \cite[Theorem 3.1]{Rawal2014}.
\end{proof}

\end{document}